\documentclass{amsart}
\usepackage{amssymb,amsmath,amsthm,latexsym}
\usepackage{amssymb,graphicx}
\newtheorem{theorem}{Theorem}

\newtheorem{conjecture}[theorem]{Conjecture}
\newtheorem{corollary}[theorem]{Corollary}

\newtheorem{definition}[theorem]{Definition}

\newtheorem{question}[theorem]{Question}
\newtheorem{proposition}[theorem]{Proposition}

\begin{document} 
\title{Algebraic properties of profinite groups}
\author{Nikolay Nikolov}
\begin{abstract}Recently there has been a lot of research and progress in profinite groups. We survey some of the new results and discuss open problems. A central theme is decompositions of finite groups into bounded products of subsets of various kinds which give rise to algebraic properties of topological groups. 
\end{abstract}
\maketitle

\section{Profinite groups}
The classification of finite simple groups has transformed the study of finite groups and this in turn has brought a wealth of results about infinite groups with various finiteness conditions. One such obvious class is the residually finite groups. Frequently questions about such groups can be reduced to asymptotic properties of their finite images and a  natural tool for studying these is the profinite groups. 

A profinite group is a compact Hausdorff topological group $G$ which is totally disconnected, i.e. any connected component of $G$ is a singleton. It is a classical result due to van Dantzig \cite{dantzig} that the last condition is equivalent to the following: any open subset of $G$ containing the identity contains an open normal subgroup. The diagonal embedding 
\[ i: G  \rightarrow \prod_{N\vartriangleleft_o G} G/N\] is a topological isomorphism and the image $i(G)$ is the inverse limit of the set $F(G)=\{ G/N \ | \ N \vartriangleleft_o G\}$ of topological finite images of $G$.

Conversely the inverse limit $\underleftarrow{\lim} \Gamma_i$ of any inverse system of finite groups $(\Gamma_i)$ is a profinite group. For details of these constructions and basic properties of profinite groups we refer the reader to \cite{RZ}, \cite{serre} or \cite{wilson}.

Profinite groups appeared first in number theory, in the first instance as a tool for studying congruences, namely the ring of $p$-adic integers $\mathbb Z_p$ and second as Galois groups of normal separable algebraic extensions. They form a part of the more general theories of compact groups and totally disconnected locally compact groups. For example the stabilizer of a vertex in the group of automorphisms of a locally finite tree is naturally a profinite group. Profinite groups feature in arithmetic geometry: the etale fundamental group $\pi(X)$ of a curve $X$ of genus at least 2 defined over a number field $k$ maps onto a $\mathrm{Gal}(k)$ and a famous conjecture of Grothendieck claims a bijection between the rational points $X(k)$ and the conjugacy classes of sections for this map, see \cite{Konig}. 

Let us mention another important example, the compact $p$-adic analytic groups. For our purposes we can define them here as the closed subgroups of $GL_n(\mathbb Z_p)$. These groups were first studied by M. Lazard \cite{lazard} as the non-archimedian equivalent of Lie groups, i.e. the compact topological groups which are analytic over $\mathbb Z_p$. This was continued by Mann and Lubotzky, in particular they obtained the following group theoretic characterization of $p$-adic groups.
\begin{theorem}[\cite{ddms}]A pro-$p$ group $G$ is $p$-adic analytic if and only if $G$ has finite rank. 
\end{theorem}
Here and below by rank we mean Pr\"{u}fer rank, i.e. the smallest integer $d$ such that any subgroup of a finite topological quotient $G/N$ of $N$ is $d$-generated.
 
The applications of profinite groups in other branches of mathematics are often via profinite completions:
Let $\Gamma$ be a residually finite group. We can define the profinite topology on $\Gamma$ by declaring the open sets to be the unions of cosets of subgroups of finite index of $\Gamma$. This makes $\Gamma$ into a Hausdorff totally disconnected group and we can define $\hat \Gamma$ to be the completion of $\Gamma$ with respect to this topology. A concrete way to construct $\hat \Gamma$ is as the closure of $i(\Gamma)$ in $\prod_{N} G/N$ where $N$ ranges over all subgroups of finite index in $G$ and $i$ is the diagonal embedding of $\Gamma$.  In this way $\Gamma$ is a subgroup of its profinite completion $\hat \Gamma$ and properties of $\Gamma$ can be deduced from those of $\hat \Gamma$. We can define the pro-$p$, pronilpotent or prosoluble completion of $\Gamma$ in a similar way using the normal subgroups $N$ of $G$ such that $G/N$ is a finite $p$-group, finite nilpotent or finite soluble group respectively.
 
A good illustration of the success of this approach is the Lubotzky's linearity condition (see \cite{ddms}, Interlude B).
\begin{theorem} A finitely generated group $\Gamma$ is linear over a field of characteristic 0 if and only if there is a chain of normal subgroups $\Gamma>\Gamma_1>\Gamma_2 \cdots$ of finite index in $\Gamma$ with $\cap_i \Gamma_i=1$ such that the inverse limit of the groups $\{\Gamma/\Gamma_i\}_i$ is a $p$-adic analytic group. Equivalently the family $\{\Gamma/\Gamma_i\}_{i}$ has bounded Pr\"{u}fer rank.
\end{theorem}
By contrast no such description is known for linear groups over fields of positive characteristic. 
Other applications of profinite groups include: the study of subgroup growth in residually finite groups \cite{LSegal}, the congruence subgroup problem \cite{alex}, the classification of $p$-groups of given coclass \cite{coclass}. 

\section{Profinite groups as profinite completions}

We explained how a finitely generated residually finite group gives rise to a profinite group, namely its profinite completion. One may ask, conversely if every profinite group is the completion of some finitely generated residually finite abstract group. A moment's thought shows that the answer is no: The $p$-adic integers $\mathbb Z_p$ are not the profinite completion of any finitely generated group. We must change the question to: Which profinite groups are completions of finitely generated groups? Let us call such a profinite group a \emph{profinite completion}. This is the same as asking: which collections of finite groups can be the images of some finitely generated residually finite group? The following example shows that we should not expect an easy answer:

\textbf{Example:} Let $G$ be a profinite group with polynomial subgroup growth which is not virtually soluble. 
A discrete (respectively profinite) group $G$ is said to have polynomial subgroup growth if there is some integer $d$ such that $G$ has at most $n^d$ (open) subgroups of index $n$ for each $n \in \mathbb N$. 

Profinite groups with polynomial subgroup growth were characterised by Segal and Shalev, see Theorem 10.3 in \cite{LSegal}. For example we can take $G$ to be a cartesian product $\prod_{i=1}^\infty PSL_2(p_i)$ for a fast increasing sequence of primes $p_i$. If $G$ is the completion of a residually finite finitely generated group $\Gamma$ then $\Gamma$ has polynomial subgroup growth as well. But then by a theorem of Mann, Lubotzky, Segal, (see \cite{LSegal} Theorem 5.1) $\Gamma$ must be a virtually soluble group of finite rank and hence so is $G$, contradiction. 

In fact a linearity criterion first proved by J. Wilson (see \cite{LSegal}, Proposition 16.4.2 (ii)) shows that $\prod_{i=1}^\infty PSL_2(p_i)$ is not a profinite completion for any infinite sequence of primes $p_i$, see \cite{dan5} for the full argument. 
\medskip

It is relaively easy to describe the nilpotent profinite groups which are profinite completions: These are the profinite groups commensurable with $\prod_{p} G(\mathbb Z_p)$ where $G$ is a unipotent algebraic group defined over $\mathbb Z$. 

There are only a few special classes of profinite groups which are known to be profinite completions, among them are iterated wreath products of finite simple groups and semisimple profinite groups. From now on we make the convention that a finite simple group means a \emph{nonabelian} finite simple group.

\begin{theorem}[\cite{dan1}]Let $S_1,S_2, \ldots, $ be a sequence of nonabelian finite simple groups, each with a transitive action as a permutation group. Form the wreath product $W_k= S_1 \wr S_2 \wr \cdots \wr S_k$ and let $G$ be the inverse limit of the groups $W_k$. Then $G$ is a profinite completion.
\end{theorem}   
In \cite{dan1} the proof assumes an extra condition on the actions of the $S_i$ but this is unnecessary, see the remark in \cite{LSegal}, page 266. 
In particular any collection of finite simple groups can be the composition factors of the finite images of some finitely generated residually finite group. 

A profinite group $G$ which is a Cartesian product of nonabelian finite simple group is called \emph{semisimple}. The next example characterises the semisimple groups which are profinite completions.

\begin{theorem}[\cite{NK}]Let $G=\prod_{i=1}^\infty S_i$ be a semisimple profinite group which is topologically finitely generated. Then 
$G$ is a profinite completion if and only if for each $n$ only finitely many of the groups $S_i$ are groups of Lie type of dimension $n$.
\end{theorem}

The dimension of $S_i$ is the dimension of the complex Lie algebra associated to the Lie type of $S_i$. \medskip

Another natural question is: Can two nonisomorphic residually finite groups have the same profinite completions? The answer is yes, even in the case of nilpotent groups, see \cite{dan3}, Chapter 11, Corollary 4. The following strong counterexample was constructed by Bridson and Grunewald \cite{BG} answering a question posed by Grothendieck:
\begin{theorem} \label{bridson} There exist two finitely presented residually finite groups $\Gamma_1, \Gamma_2$ with an injection $i: \Gamma_1 \rightarrow \Gamma_2$ such that 
the induced map $\hat i: \hat \Gamma_1 \rightarrow \hat \Gamma_2$ is an isomorphism but $\Gamma_1$ and $\Gamma_2$ are not isomorphic. 
\end{theorem} 
In \cite{longreid} the authors show that counterexamples to Grothendieck's question as in Theorem \ref{bridson} cannot exist within arithmetic hyperbolic 3-manifold groups and pose the following
\begin{question} Suppose that $M_1$ and $M_2$ are geometric 3-manifolds with infinite fundamental group for which the
profinite completions $\widehat{\pi_1(M_1)}$ and $\widehat{\pi_1(M_2)}$ are isomorphic. Are $M_1$ and $M_2$ homeomorphic?
\end{question}

It turns out that there exist uncountably many pairwise non-isomorphic finitely generated residually finite groups with the same profinite completion, see \cite{laci}. On the other hand it is a theorem of Grunewald, Pickel and Segal \cite{GPS} that there are only finitely many polycyclic groups with a given profinite completion. Similar result has been proved by Aka \cite{aka} for simple arithmetic groups with the congruence subgroup property. However when we turn to more general classes of groups even the following is open, see \cite{GruZal}:
\begin{question}[Remeslennikov]\label{hat}
A finitely generated residually finite group $\Gamma$ has the same profinite completion as the free group on two generators $F_2$. Must $\Gamma$ be isomorphic to $F_2$?
\end{question}
If $\Gamma$ is two generated then $\Gamma$ must indeed be isomorphic to $F_2$: Any surjection $\pi: F_2 \rightarrow \Gamma$ gves rise to a surjection $\hat \pi: \hat F_2 \rightarrow \hat \Gamma \simeq \hat F_2$ of their profinite completions. Now a finitely generated profinite group $G$ is non-Hopfian, i.e. any surjection from $G$ to $G$ is an isomorphisms. Therefore $\hat \pi$ is an isomorphism and hence so is $\pi$. So Question \ref{hat} is really about groups $\Gamma$ which need more than 2 generators. 

We remark that the analogous question for pro-$p$ completions has a negative answer: There exist groups $\Gamma$ which are not free but have the same set of nilpotent images as a free group. Such groups are called \emph{parafree}, for examples see \cite{baum}.
\medskip

When we consider profinite completions of finitely presented groups there are a few more restrictions. First let us make two definitions. A group $G$ is \emph{large} if it has a subgroup of finite index $H$ which maps homomorphically onto a nonabelian free group. Similarly $G$ is said to be $p$-large for a prime $p$ if $H$ above can be taken to be a subnormal subgroup with $[G:H]$ a power of $p$. Clearly a large group $G$ will have a lot of finite images, in particular any finite group will appear as an image of a finite index subgroup of $G$. 

For a prime $p$ a chain of subgroups $(G_i)$ is called an abelian $p$-series of \emph{rapid descent} if $G_i > G_{i+1} \geq G_i^p[G_i,G_i]$ for all $i \in \mathbb N$ and \[ \liminf_{i \rightarrow \infty} \frac{\dim_{\mathbb F_p} G_i/G_{i+1}}{[G:G_i]} >0.\]
The following theorem has been proved by M. Lackenby in \cite{detect}.
 
\begin{theorem}\label{marc} A finitely presented group $G$ has an abelian $p$-series of rapid descent if and only if $G$ is $p$-large.
\end{theorem}

As a corollary if two finitely presented groups have isomorphic profinite completions then one of them is large if and only if the other is large.
 
\section{Rigidity of profinite groups}
Let $G$ be a profinite group. What happens if we take the profinite completion of $G$ itself? Since the open subgroups of $G$ have finite index the profinite topology on $G$ contains the original topology but could it be strictly stronger? In \cite{serre} Serre asked this question in the form: Assuming that $G$ is (topologically) finitely generated is it true that every finite index subgroup of $G$ open? Here and below we say that a profinite group is finitely generated if it contains a dense finitely generated subgroup.

Let us call a profinite group $G$ \emph{rigid} (or strongly complete) if every subgroup of finite index is open in $G$. 
It is easy to see that there exist non-rigid groups which are not finitely generated. For example if $G=C_p^{\aleph_0}$ then $G$ has $2^{2^{\aleph_0}}$ subgroups of index $p$ of which only $\aleph_0$ are open. Even worse: the same abstract group can support two inequivalent topologies as a profinite group: The group 
$G=\prod_{i=1}^\infty C_{p^i}$ is abstractly isomorphic to $G_1=G \times \mathbb Z_p$ but there is no continuous isomorphism between $G$ and $G_1$. More generally in \cite{jonny} J. Kiehlmann has classified all countably based abelian pro-$p$ groups up to continuous and up to abstract isomorphisms. 

When $G$ is a finitely generated pro-$p$ group then Serre himself proved in \cite{serre} that $G$ is rigid. Increasing general cases were proved by Anderson \cite{anderson}, Hartley \cite{hartley}, Segal \cite{dan4} until in 2003 the author and D. Segal answered Serre's question in the positive:

\begin{theorem}[\cite{NS3}] \label{rigid} Every finitely generated profinite group is rigid. 
\end{theorem}

In particular all homomorphisms between finitely generated profinite groups are automatically continuous and each such group is its own profinite completion.

\section{Word width in finite and profinite groups}
The proof of Theorem \ref{rigid} relies on \emph{bounded width} of commutators and other words in finite groups.
  
It is a standart result observed by Brian Hartley that algebraic properties of a profinite group $G$ are equivalent to asymptotic properties of the collection $F(G)$ of finite continuous images of $G$. For example $G$ is topologically finitely generated if and only if there is some $d$ such that any member of $F(G)$ is generated by $d$ elements. For a subset $X$ of $G$ let us write \[ X^{*n}= \{ x_1 \cdots x_n \ | \ x_i \in G\}\] The following proposition is crucially used when proving that some subgroups of $G$ are closed.

\begin{proposition} \label{asy} Let $X=X^{-1}$ be a closed subset of a profinite group $G$. The group $H=\langle X\rangle$ is closed in $G$ if and only if there is some $n \in \mathbb N$ such that $H=X^{*n}$. In turn this condition is equivalent to $\bar X ^{*n}=\bar H$ for the images $\bar X, \bar H$ of $X$ and $H$ in every $\bar G \in F(G)$.
\end{proposition}

A typical application is where $X$ is the set $G_q:=\{x^q\ | \ x\in G\}$ of $q$-powers in $G$ or the set of commutators $[x,y]$ in $G$. We can summarize the main results of \cite{NS3} and \cite{NS2} as follows.
\begin{theorem}\label{powers} Let $d, q \in \mathbb N$ and let $\Gamma$ be a finite $d$-denerated group. There exist functions $f_1(d,q)$ and $f_2(d)$ such that

1. Every element of $\Gamma^q$ is a product of at most $f_1$ powers $x^q$. i.e. $\Gamma^q=\Gamma_q^{*f_1}$.

2. For a normal subgroup $N$ of $\Gamma$ the group $[N,\Gamma]$ generated by the set $Y=\{ [n,g] \ | \ n \in N, g \in \Gamma\}$ is equal to $Y^{*f_2}$.
\end{theorem}
The proof of this theorem has recently been streamlined in \cite{NS} and we shall indicate some of the main ingredients in section \ref{sec:images} below. 
Proposition \ref{asy} now gives
\begin{corollary} \label{cor1} Let $G$ be a finitely generated profinite group. For any integer $q$ and any closed normal subgroup $N$ of $G$ the algebraically defined subgroups $G^q$ and $[N,G]$ are closed in $G$. In particular the terms of the lower central series $\{\gamma_i(G)\}_{i=1}^\infty$ of $G$ are closed.
\end{corollary}
Corollary \ref{cor1} implies that $G/G^q$ is a profinite group and therefore is an inverse limit of finite $d$-generated groups of exponent $q$. By the solution of the restricted Burnside problem by Zelmanov \cite{zelmanov} the sizes of these finite groups are bounded, i.e. $G/G^q$ is finite and therefore $G^q$ is open. 
Theorem \ref{rigid} is now an easy consequence: 

Let $H$ be a subgroup of index $q$ in a finitely generated profinite group $G$. We want to prove that $H$ is open. Without loss of generality we may assume that $H$ is normal in $G$ and hence $H \geq G^q$, which is open in $G$ by the above argument. Hence $H$ is open in $G$.

In fact the use of Zelmanov's theorem above was not strictly speaking necessary and it can be avoided if one just wants to prove Theorem \ref{rigid}. This was done
in \cite{NS3} (which was published before it was known that $G^q$ is closed), and a simplified argument can be found in \cite{NS}, section 5.1, using Theorem \ref{fin} below.

One might ask whether Corollary \ref{cor1} extends to other algebraically defined subgroups of $G$, for example verbal subgroups. Let us make this more precise. Let $w=w(x_1, \ldots x_k)$ be an element of the free group on $x_1, \ldots ,x_k$, we shall refer to $w$ as a word in $x_i$. The set \[ G_w:=\{w(g_1, \ldots, g_k)^{\pm 1} \ | \ g_i \in G \} \] is called the set of values of $w$ in $G$. It is clearly a closed subset of $G$ when $G$ is a profinite group. The verbal subgroup $w(G)$ is defined as $w(G):=\langle G_w\rangle$. The word $w$ is said to have width at most $n$ in a group $G$ if $w(G)=G_w^{*n}$. Then the verbal subgroup $w(G)$ is closed in $G$ if and only if there is some $n \in \mathbb N$ such that $w(\bar G)= \bar G_w^{*n}$ for all finite continuous images $\bar G$ in $F(G)$. If this is so we shall say that $w$ has bounded width in the family $F(G)$.

\begin{question} \label{words} For which words $w$ it is true that $w(G)$ is closed in all finitely generated profinite groups $G$?
\end{question} 

This is equivalent to $w$ having bounded width in all finite $d$-generated groups. Theorem \ref{powers} provides us with many such words:

\begin{corollary} Let $d \in \mathbb N$. Suppose that $w$ is either a non-commutator word (i.e. $w\not \in F'$) or the word $w=[x_1,x_2, \ldots, x_m]$. Then $w$ has bounded width in all finite $d$-generated groups. 
\end{corollary}

Could the above corollary hold for all words? This is not true: Romankov \cite{romankov} gave an example of $3$-generated pro-$p$ group $G$ such that the second derived group $G''$ is not closed, i.e. the width of the word $[[x_1,x_2],[x_3,x_4]]$ is unbounded in $d$-generated finite $p$-groups. The most definitive describition so far has been achieved by A. Jaikin-Zapirain \cite{andrei} who has answered Question \ref{words} for pro-$p$ groups.

\begin{theorem} Let $w \not = 1$ be a word in a free group $F$. The following are equivalent:

1. $w(G)$ is closed in all finitely generated pro-$p$ groups $G$.

2. $w \not \in F''(F')^p$. 
\end{theorem}   

Let us call a word $w$ a $J$-word if $w \not \in F''(F')^p$ for any prime $p$. It is clear that $w$ must be a $J$-word for $w(G)$ to be closed in each finitely generated profinite group $G$. Is the converse true?  No counterexample is known to this question. D. Segal \cite{dan} has proved that if $w$ is a $J$-word then $w(G)$ is closed in all finitely generated soluble profinite groups $G$. Natural candidates to try next are the Engel words $[x_1,x_2,\ldots,x_2]$. 
\medskip

Could Question \ref{words} have a positive answer for all words in some restricted class of groups? It turns out that this is the case if we consider only $p$-adic analytic groups:
\begin{theorem} \cite{andrei} Let $G$ be a compact $p$-adic analytic group and let $w$ be any word. Then $w(G)$ is closed.
\end{theorem}
L. Pyber has asked if this holds in greater generality for finitely generated adelic groups, i.e. closed subgroups of $\prod_p GL_n(\mathbb Z_p)$. Dan Segal \cite{dan2} has answered this question in the positive for adelic groups which have the property FAb, i.e. every open subgroup has finite abelianization. A related question is:

\begin{question} Let $w$ be a word and $r \in \mathbb N$. Is there $f=f(w,r)$ such that for any $p$-adic analytic group $G$ of rank $r$ the width of $w$ in $G$ is at most $f(w,r)$?  
\end{question}
The other extreme to $p$-groups is the family of finite simple groups. A finite simple group is generated by $2$ elements, so Thereom \ref{powers} implies that the word width of any power word $x_1^q$ or $[x_1,x_2]$ in the finite simple groups is bounded (a result proved in \cite{MZ} and \cite{SW} earlier and in fact an essential ingredient in the proof of Theorem \ref{powers}). Much better bounds are now known for \emph{all} words. 
\begin{theorem}[\cite{LST}] \label{aner} Let $w$ be a word and $S$ be a finite simple group. If the size of $S$ is big enough as a function of $w$ alone, then $S=S_w S_w$, i.e. $w$ has width 2 in $S$.
\end{theorem}

It follows that $w(G)$ is closed in all semisimple profinite groups $G$.

The obvious example $w=x_1^q$ with $q$ dividing $|S|$ shows that 2 is the best possible bound in Theorem \ref{aner} in general. This is the essentially the only example we know where the word map is not surjective on finite simple groups. The following has been proved by Liebeck, Shalev, O'Brien and Tiep \cite{ore}.

\begin{theorem}[Ore Conjecture] If $S$ is a finite simple group then every element of $S$ is a commutator.
\end{theorem}

A. Shalev (private communication) has suggested the following problem: Let $w$ be a word which is not a proper power. If $S$ is a finite simple group of Lie type of large enough rank or a large alternating group then $S=S_w$.

\subsection{Word width in discrete groups}
Moving away from profinite groups let us say a few words about words in abstract groups. One of the first people to investigate word width was Philip Hall, who formulated a series of conjectures about verbal subgroups, see \cite{robinson}, Chapter 4.2. Hall's definition of elliptic words is the same as our definition of words of finite width we have adopted here. Hall's student P. Stroud \cite{stroud} proved that any word has finite width in an abelian-by-nilpotent group. Independently Romankov \cite{romankov} proved the same result for virtually polycyclic groups.
D. Segal \cite{dan} recently extended this and proved that the word width in virtually soluble minimax groups is always finite. 

One might expect that most words will have infinite width in a free group and this is indeed correct: A. Rhemtulla \cite{rhem} showed that with the exception of trivial cases a word has infinite width in free products. This has recently been generalized to hyperbolic groups \cite{hyper}. Some open problems remain:

\begin{question} Let $\Gamma$ be a centre by metabelian group (i.e. $\Gamma/Z(\Gamma)$ is metabelian). Is it true that some word $w$ has infinite width in $\Gamma$? 
\end{question} 
\begin{question} Let $\Gamma$ be a finitely generated soluble group. Is it true that the commutator word $[x_1,x_2]$ has finite width in $\Gamma$?
\end{question}
This has been answered affirmatively by Rhemtulla \cite{rhem2} when $\Gamma$ has derived length at most 3. It is open for the free soluble group on two generators of derived length 4.

\section{Fibres of word maps}
Let $\Gamma$ be a finite or infinite group. We can consider a word $w \in F_k$ as a map $(g_1,\ldots, g_k) \rightarrow w(g_1, \ldots, g_k)$ from $\Gamma^{(k)} \rightarrow \Gamma$ and investigate the properties of this map. 

For an element $g \in \Gamma$ let us write 
\[ P_\Gamma(w,g)= \frac{|\{ \mathbf x \in \Gamma^{(k)} \ | \ w(\mathbf x)=g \}|}{|\Gamma|^k} \] 
This is the probability of satisfying $w(\mathbf x)=g$ in $\Gamma$ for a random $k$-tuple $x \in \Gamma^{(k)}$.

A word $w$ is said to be \emph{measure preserving} (in finite groups) if $P_\Gamma(w,g)=|\Gamma|^{-1}$ for all finite groups $\Gamma$ and all $g \in \Gamma$. An example is any primitive word, i.e. an element of a free basis of the free group $F_k$. In fact a word $w$ is measure preserving if and only if $w$ is a primitive element in the free profinite group $\hat F_k$, see \cite{pud}, Section 6.

T. Gelander has conjectured that conversely any measure preserving word must be primitive in $F_k$. D. Puder \cite{pud} proved this in the case of $F_2$. Very recently this has been extended by D. Puder and O. Parzanchevski to all words.
\begin{theorem}[\cite{pp}] Let $w \in F_k$. Then $w$ is measure preserving in all finite groups if and only if $w$ is primitive.
\end{theorem}

Let us mention another intriguing open problem concerning the fibres of the word map.
If $\Gamma$ is abelian then the word map is a homomorphism and $P_\Gamma(w,e)= |w(\Gamma)|^{-1} \geq |\Gamma|^{-1}$, in particular $P_\Gamma(w, e)$ is bounded away from $0$ for any word $w$. A. Amit (unpublished) asked whether this property characterizes all soluble groups. The combined results in \cite{Ab1} and \cite{ns1} answer this affirmatively (without using the classification):

\begin{theorem} Let $\Gamma$ be a finite group. Then $\Gamma$ is soluble if and only if there exists $\epsilon>0$ such that for any word $w\in F_k$ we have $P_\Gamma(w,e) \geq \epsilon$. 
\end{theorem}
 Amit raised the following
 \begin{question} Is it true that when $\Gamma$ is a finite nilpotent group and $w\in F_k$ is a word then $P_\Gamma(w,e) \geq |\Gamma|^{-1}$?
 \end{question}
M. Levy \cite{Levy} has answered this in the positive when $\Gamma$ has nilpotency class 2. The general case remains open.

When we consider $P_\Gamma(w,g)$ for a finite simple group $\Gamma$, Larsen and Shalev \cite{larsenshalev}have proved the following with applications to subgroup growth and representation varieties. 
\begin{theorem} Given a word $w$ there is $\epsilon= \epsilon(w)>0$ and $N=N(w)$ such that for any finite simple group $S$ of size at least $N$ and any $ g \in S$ we have  $P_S(w,g)< |S|^{-\epsilon}$.
\end{theorem}
It turns out that some words are \emph{almost measure preserving} in the finite simple groups. In order to define what this means,  let us write for $X \subset G$ 
$P_G(w,X)= \sum_{g\in X} P_G(w,g)$.
 
\begin{definition} A word $w \in F_k$ is said to be almost measure preserving in a family $\mathcal A$ of finite groups if for any $\epsilon>0$ there is $N=N(\epsilon)$ such that for any $G \in \mathcal A$ with $|G|>N$ and any subset $X \subset G$ we have
\[  |P_G(w,X) - \frac{|X|}{|G|}| < \epsilon. \]
\end{definition}
In \cite{shelly} S. Garion and A. Shalev show that the commutator word $[x_1,x_2]$ is almost measure preserving in all finite simple groups. This has diverse applications, in particular to Theorem \ref{T1} below. Similar result for the words $x_1^nx_2^m$ have been proved by M. Larsen and A. Shalev (work in preparation).
\section{Profinite groups of type IF}
A finitely generated profinite group has only finitely many open subgroups of any given index. Let is call a profinite group with the latter property a group of type IF. Similarly let us define a profinite group to be of type AF if it has finitely many subgroup of any given finite index (which  may or may not be open in $G$). Clearly AF implies IF and IF is equivalent to AF for rigid profinite groups. 

A profinite group of type AF may not be finitely generated as the example $H=\prod_{n \geq 5} A_n^{(n!)^n}$ shows: It is clear that $H$ has type IF while on the other hand from Theorem \ref{aner} it follows that $H^q$ is open in $H$ for any $q$ which shows that $H$ is rigid. 

It is not too hard to show that a rigid profinite group $G$ must have AF (and so also IF), a result first proved in \cite{Pit}. Indeed otherwise $G$ will maps topologically onto a Cartesian product $L=\prod_{i=1}^\infty S_i$ of isomorphic finite simple groups $S_i \simeq S$ (maybe abelian). Choose a non-principal ultrafilter $\mathcal U$ on $\mathbb N$ and consider the ultraproduct $L/K$ where \[ K= \{ (g_i) \in L  \ | \ \textrm{ the set } \{ i \ | \ g_i=e \} \in \mathcal U\}. \] 
Then $L/K \simeq S$ and so $K$ is a normal subgroup of finite index in $L$ which is not open. Thus $L$ is not rigid and neither is $G$. 
In fact \cite{SmWil} proves that a profinite group $G$ is rigid if and only if it has type AF. One is led therefore to ask whether Theorem \ref{rigid} can be generalized to all IF-groups, i.e. could it be that all IF-groups are rigid. However this is not true:
\begin{proposition} There exists a group of type IF which is not rigid.
\end{proposition}

The idea of the proof is to use square width: Define $sw(G)$ to be the smallest integer $k$ (if it exists) such that any element of $G^2$ is a product of $k$ squares and set $sw(G)=\infty$ otherwise. We shall find a sequence of finite groups $\{\Gamma_n\}_{n=5}^\infty$ such that:

1. $\Gamma_n=\Gamma_n'$

2. $\Gamma_n$ has a unique maximal normal subgroup $K_n$ such that $\Gamma_n/K_n \simeq A_n$.

3. $sw(\Gamma_n)\geq n$. 

If we then take $G=\prod_{n=5}^\infty \Gamma_n$ then condition $2$ implies that $G$ has type IF. On the other hand condition 3 and Proposition \ref{asy} imply that $G^2 \not = G$ hence $G$ has a subgroup of index $2$ which cannot be open by condition 1.

Let $Q_n$ be a perfect finite group with $sw(Q_n)>2n^2$. Such a group exists for any $n$ by \cite{Ho}, Lemma 2.2. Put $\Gamma_n= Q_n \wr A_n= Q_n^{(n)} \rtimes A_n$. It is clear that $\Gamma_n$ satisfies 1 and 2. As for the third condition I claim that $cw(\Gamma_n)\geq  cw(Q_n)/2n \geq n$. 

To prove this let $s= cw(\Gamma_n)$ and express the element $\mathbf{g}=(x,1,1, \cdots, 1) \in Q_n^{(n)}$ as a product of $s$ squares in $\Gamma_n$:
\[ \mathbf{g} = \prod_{i=1}^s (\mathbf b_i \pi_i )^2, \quad \mathbf{b}_i =(b_i(1), \ldots ,b_i(n))\in Q_n^{(n)}, \quad \pi_i \in A_n.\]
By collecting the elements of the base group to the left we reach the equation $\mathbf g= \prod_{i=1}^s \mathbf{b_i}^{\alpha_i} \mathbf{b_i}^{\beta_i}$ where $\alpha_i,\beta_i$ are some permutations from $A_n$. Now if we multiply together all the coordinates of the elements on both sides of this equation we reach
the equation $x=U$ where $U$ is product of elements $b_i(j) \in Q_n$ ($i=1, \ldots, s, \ j=1, \ldots, n$) in some order, each appearing exactly twice.

Now we only need to observe that by the proposition below $U$ is a product of at most $2sn$ squares. Since $x \in Q_n$ was arbitrary this implies that $sw(Q_n) \leq 2ns = 2n \cdot sw(\Gamma_n)$ which was what we were after.
\begin{proposition} Let $\Gamma$ be a group and let $U$ be a product of length $2m$ of elements $c_1, \ldots ,c_m$ in some order, each appearing exactly twice. Then $U$ is a product of $2m-1$ squares.
\end{proposition}
\textbf{Proof:} Suppose that $U=c_rV_1c_rV_2$, where $V_1V_2$ does not involve $c_r$. Then 
\[ U=(c_rV_1)^2 V_1^{-2} V_1V_2\]
and by induction we may assume that $V_1V_2$ is a product of $2m-3$ squares. $\square$
\section{Presentations and cohomology}
Let $\Gamma$ be a finite group with generating set of minimal size equal to $d$. Then $\Gamma=F_d/N$ is a quotient of the free group $F_d$ on the set $X=\{x_1, \ldots, x_d\}$ by a normal subgroup $N$. By a minimal presentation for $\Gamma$ we mean a presentation $\langle X | \ R \rangle$ where $R$ is a set of relators with $|R|$ as small as possible and we define $r(\Gamma)=|R|$. 

At the same time $\Gamma$ is a quotient of the free profinite group on $X$, $\hat F_d$. By a profinite presentation of $\Gamma$ we mean a pair $\langle X | R_1\rangle$ where $R_1 \subset \hat F_d$ such that $\Gamma \simeq \hat F_d/ U$, with $U:=\overline{\langle R_1^{\hat F_d}\rangle}$, the closed normal subgroup of $\hat F_d$ generated by $R_1$. Again a minimal profinite presentation is one where $|R_1|$ is as small as possible and we set $\hat r (\Gamma)=|R_1|$. We can view every abstract presentation of $\Gamma$ as profinite presentation and this shows that $\hat r(\Gamma) \leq r(\Gamma)$. It is a well-known problem whether this is always an equality.
\begin{question} Is there a finite group $\Gamma$ with $\hat r(\Gamma) < r(\Gamma)$?
\end{question}
The motivation for this question is that unlike $r(\Gamma)$ the quantity $\hat r(\Gamma)$ can in theory be computed from the representation theory of $\Gamma$, a classic result of Gruenberg, for the formula see Proposition 16.4.7 of \cite{LSegal}. In particular for a finite $p$-group $P$ $\hat r (P)=H^2(P, \mathbb F_p)$.

A case of interest is when $\Gamma$ is a finite simple group: In \cite{GKKL} the authors show that $\hat r(\Gamma) \leq 18$ and $\hat r (A_n) \leq 4$. As for abstract presentation of finite simple groups in another paper \cite{GKKL2} the same authors prove that $r(\Gamma) \leq 80$ with $r(A_n) \leq 8$.
\medskip

Let us return to Serre's result that any finite index subgroup in a finitely generated pro-$p$ group is open. It can be restated in the following way: any homomorphism from $G$ to $(\mathbb F_p,+)$ is continuous.

For a profinite group $G$ and a finite topological $G$-module $M$ (which just means that $C_G(M)$ is open in $G$) define $H^n_c(G,M)$ to be usual cohomology group defined as the quotient $Z^n_c(G,M)/B^n_c(G,M)$ of the groups of continuous $n$-cocycles and $n$-coboundaries from $G$ to $M$. Similarly, let $H_a^n(G,M)$ be the analogue defined without requiring continuity of the cocycles or coboundaries. It is easy to see that $H_c^n(G,M)$ embeds in $H_a^n(G,M)$. For example when $M$ is a trivial module $H_c^1(G,M)$ is the additive group of continuous homomorphisms from $G$ to $(M,+)$, while $H_a^1(G,M)$ is the group of all abstract homomorphisms. 

Serre's result can now be restated as $H_c^1(G,F_p)=H_a^1(G,F_p)$ provided $H_c^1(G,F_p)$ is finite. It is natural to ask if the higher dimensional analogue of this holds. 
\begin{question} \label{homology}Let $G$ be a finitely presented pro-$p$ group (i.e. $H^1_c(G,F_p)$ and $H_c^2(G,F_2)$ are finite). Is it true that $H_a^2(G,F_p)=H_c^2(G,F_p)$?
\end{question}

It is not too hard to see that the condition thet $G$ is finitely presented is necessary, see \cite{alcober}. A case of special interest is when $G$ is a $p$-adic analytic group, which has been answered affirmatively for Chevalley and soluble $p$-adic groups by \cite{sury}. Sury's result concerns non-compact Chevalley groups over $\mathbb Q_p$ and his methods have been extended to cover the compact Chevalley group $SL_2(\mathbb Z_p)$ as well by Barnea, Jaikin-Zapirain and Klopsch (work in preparation). Another of their results is that a positive answer to Question \ref{homology} for free pro-$p$ groups implies a positive answer in general for all finitely presented pro-$p$ groups.

The groups $H^2(G,\mathbb F_p)$ (abstract or continuous) parametrize the equivalence classes of (abstract or continuous) extensions of $F_p$ by $G$. Thus the following are equivalent for a pro-$p$ group $G$: \medskip

1. $H_a^2(G,F_p)=H_c^2(G,F_p)$.

2. Any central extension of $C_p$ by $G$ (i.e. a group $K$ with a normal central subroup $C$ of order $p$ such that $K/C \simeq G$) is isomorphic to a pro-$p$ group. 

3. Any central extension of $C_p$ by $G$ is residually finite.
\medskip

Using this now it is immedate that $H_a^2(G,F_p)\not =H_c^2(G,F_p)$ when $G$ is not finitely presented: Let $G=F/N$ where $F$ is a finitely generated free pro-$p$ group and $N$ is a closed normal subgroup. If $G$ is not finitely presented then $N/[N,F]N^p$ is infinite elementary abelian pro-$p$ group and hence it has a subgroup $N>K>[N,F]N^p$ of index $p$ in $N$ which is not closed. This $F/K$ is an extension of $C_p=N/K$ by $G=F/N$ which is not residually finite. 

Similarly one can prove the following:

\begin{proposition} \label{ext} Let $G$ be a pro-$p$ group which is an extension of finitely generated pro-$p$ group $N$ by a finitely generated profinite group $H=G/N$. Suppose that $H_a^2(N,F_p)=H_c^2(N,F_p)$ and $H_a^2(U,F_p)=H_c^2(U,F_p)$, for any open subgroup $U$ of $H$. Then $H_a^2(G,F_p)=H_c^2(G,F_p)$.
\end{proposition}

\textbf{Proof:} First note that the condition on $H$ implies that any extension $E$ of a finite $p$-group $P$ by $H$ is residually finite and so topological. This is proved by induction on $|P|$ the case when $P=C_p$ being part of the assumptions. For the induction step let $Z$ be the centre of $P$. Then $H=E/P$ acts on the finite group $Z$, so by replacing $H$ by an open subgroup (and $E$ by a finite index subgroup) we may assume that $Z \leq Z(E)$. Let $C$ be a subgroup of order $p$ in $Z$, then by considering $E/C$ and the induction hypothesis we may assume that $E/C$ is residually finite and in particular $E$ has a finite index subgroup $E_1$ with $E_1 \cap P= C$. Then $E_1/C \simeq E_1P/P$ is isomorphic to an open subgroup of $H$, so by assumption $E_1$ is a redisually finite extension, so there is a finite index  subgroup $K$ of $E_1$ with $C \cap K =1$. This implies that $E$ itself is residually finite. 

Now we can easily finish the proof: Let $C_p \vartriangleleft J \vartriangleleft K$ be an extension of $C_p \leq Z(K)$ by $G \simeq K/C_p$ where $K/J \simeq H$ and $J/C_p \simeq N$. Since $G$ is residually finite it is enough to find a subgroup $M$ of finite index in $K$ with $C_p \cap M=1$. From the assumption on $N$ we know that $J$ is residually finite and finitely generated topologically, thus we may find a subgroup $S \leq J$ with $[J:S] <\infty$, such that $S \vartriangleleft K$ and $S \cap C_p=1$. Consider now $K/S$ which is an extension of the finite $p$-group $J/S$ by $H\simeq K/J$ and so by the argument above there is some finite index subgroup $M/S$ in $K/S$ with $1=M/S \cap (C_pS)/S =(M \cap C_pS)/S$, this $C_p \cap M=1$.
$\square$

Proposition \ref{ext} shows that Question \ref{homology} for $p$-adic analytic groups reduces to the case when $G$ is a simple $p$-adic analytic group i.e. whose Lie algebra is simple. A first test case will be $G=SD^1(\Delta_p)$, the group of norm 1 elements in the quaternion division algebra over $\mathbb Z_p$. For definition of this group see Exercise 9.3 in Chapter I of \cite{KNV}.

\section{Strange images of profinite groups} \label{sec:images}
In this section we present some recent results of D. Segal and the author in \cite{NS}.
Let $G$ be a finitely generated profinite group. We can interpret Theorem \ref{rigid} as saying that every finite quotient of $G$ is topological. What can we say about other quotients? Suppose first that $G/N$ is a residually finite quotient, then it must be a profinite group. Indeed $N$ is an intersection of finite index subgroups each of which is open and hence closed in $G$. Therefore $N$ is a closed subgroup of $G$ i.e. the induced topology on $G/N$ makes it into a profinite group. As a consequence $G$ cannot have a countably infinite residually finite quotient. 

Could it be that $G$ has a countably infinite quotient? The answer is perhaps surprisingly yes. For example let $G= \prod_{p \in P}\mathbb F_p$ (product over the set $P$ of all primes $p$). Take a non-principal ultrafilter on $P$  and let $G/K$ be the ultraproduct. Then $G/K$ is an field of characteristic 0, which maps onto $\mathbb Q$ as an additive group. Therefore $(\mathbb Q,+)$ is an image of $G$. Similar argument shows that $\mathbb Z_p$ maps onto $\mathbb Q$ (use that $\mathbb Q_p$ is a vector space over $\mathbb Q$). More generally any infinite abelian profinite group has a countable infinite image. 

Let us then put a further restriction on the image: Could a profinite group have a finitely generated infinite image? The answer is no, even more generally for compact Hausdorff groups:

\begin{theorem}\cite{NS}\label{images}
Let $G$ be a compact group and $N$ a normal subgroup of (the underlying
abstract group) $G$ such that $G/N$ is finitely generated. Then $G/N$ is finite.
\end{theorem}

We can in addition provide some information about possible countable images:

\begin{theorem}
Let $G$ be a finitely generated profinite group. Let $N$ be a normal subgroup of (the
underlying abstract group) $G$. If $G/N$ is countably infinite then $G/N$ has
an infinite virtually-abelian quotient.
\end{theorem}

The last condition implies that $G$ has an open subgroup $K$ such that $K/K'$ is infinite.
This is also easily seen to be sufficient for the existence of countable infinite quotients:

\begin{corollary} Let $G$ be a finitely generated profinite group. Then $G$ has a countably infinite 
quotient if and only if some open subgroup of $G$ has infinite abelianization. 
\end{corollary} 

Note that if $G/N$ is countable then $\bar N$ is open in $G$. We say that a normal subgroup $N$ is \emph{virtually dense} in $G$ if $[G:N]$ is infinite and the closure of $N$ is open in $G$. For example, when $G$ is abelian or semisimple then $G$ has a dense normal subgroup. In the latter case we can take $N= \oplus_i S_i$ in $G= \prod_i S_i$. In fact the existence of virtually dense normal subgroups is described by these examples.

\begin{corollary} A finitely generated profinite group has a virtually dense normal subgroup if and only if it has an open normal subgroup $H$ which has an infinite abelian quotient or an infinite semisimple quotient. 
\end{corollary}
The key to proving the above results is the following result about finite groups proved in \cite{NS}.
For a finite group $\Gamma$ we denote the derived group by $\Gamma'$ and write $\Gamma_{0}$ for the
intersection of the centralizers of the non-abelian simple chief factors of
$\Gamma$.

\begin{theorem}\cite{NS} \label{fin}
Let $\Gamma$ be a finite group, $H\leq \Gamma_{0}$ a normal subgroup of $\Gamma$, and $\{y_{1}%
,\ldots,y_{r}\}$ a symmetric subset of $\Gamma$. If $H\left\langle y_{1}%
,\ldots,y_{r}\right\rangle =\Gamma^{\prime}\left\langle y_{1},\ldots,y_{r}%
\right\rangle =\Gamma$ then%
\[
\lbrack H,\Gamma]=\left(
{\displaystyle\prod\limits_{i=1}^{r}}
[H,y_{i}]\right)  ^{\ast f}%
\]
where $f=f(r,\mathrm{d}(\Gamma))=O(r^{6}\mathrm{d}(\Gamma)^{6})$.

If $\{y_{1},\ldots,y_{r}\}$ actually generates $\Gamma$ then the above equality holds
for every $H$ (not necessarily inside $\Gamma_{0}$).
\end{theorem}

Now $\Gamma/\Gamma_{0}$ is semisimple-by-(soluble of derived length
$\leq3$), while $\Gamma/\Gamma^{\prime}$ is abelian: so the theorem reduces certain
problems to the case of semisimple groups and abelian groups.

As an application let us indicate how to deduce Theorem \ref{powers} from Theorem \ref{fin}. Part (2) of Theorem \ref{powers} follows just by setting $H=N$ and $y_1,\ldots y_d$ to be a symmetric generating set of $\Gamma$.
Part (1) of Theorem \ref{powers} needs a little more work. For a finite $d$-generated group $\Delta$ we set $\Gamma:= \Delta^q$. The index of $\Gamma$ in $\Delta$ is bounded in terms of $d$ and $q$ and so $\Gamma$ can be generated by some $f_1(d,q)$ elements. Moreover we can choose a symmetric set $y_1, \ldots, y_r$ with $r$ bounded in terms of $d,q$ such that each $y_i=x_i^q$ for some $x_i \in \Delta$ and $y_1, \ldots ,y_r$ generates $\Gamma$ modulo $\Gamma'$ and modulo $\Gamma_0$. (Using a version of Theorem \ref{powers} for finite semisimple and finite soluble groups which were proved earlier) 
Then Theorem \ref{fin} with $H=\Gamma_0$ implies that every element of $[\Gamma_0,\Gamma]$ is a product of boundedly many copies of $[\Gamma_0,y_i]$. Since $[x,y_i]= y_i^{-x}y_i$ is a product of two $q$-th powers it follows that $[\Gamma_0, \Gamma]$ is a product of boundedly many $q$-th powers. It remains to deal with $\Gamma/[\Gamma_0,\Gamma]$ which is easy. The full details can be found in \cite{NS} Theorem ??.

An easy consequence of Theorem \ref{fin} is
\begin{corollary} \label{im} Let $G$ be a finitely generated profinite group and $N$ a normal subgroup of $G$. If $NG'=NG_0=G$ then $N=G$.
Here $G_0$ is the intersection of all open normal subgroups $M$ such that $G/M$ is almost simple.
\end{corollary}

Indeed, we may choose elements $y_1, \ldots, y_{2d}$ from $N$ such that $\overline{ \langle y_1, \ldots, y_{2d} \rangle }$ projects onto $G/G'$ and $G/G_0$. The profinite analogue of Theorem \ref{fin} now implies that for the integer $f$ defined there we have 
\[ [G_0,G] = (\prod_{i} [G_0,y_i])^{*f} \]
In particular $[G_0,G] \leq N$ since $y_i \in N$. Therefore \[ G'=[G,NG_0] \leq [G,N][G,G_0] \leq N\] hence $G=NG' =N$. \medskip

Corollary \ref{im} shows that if $G$ has a `strange' normal subgroup, then either
$G/G^{\prime}$ or $G/G_{0}$ has one: and these groups are not so hard to understand.
\medskip
Let us illustrate this by proving Theorem \ref{images} at least for profinite groups $G$. The general case is obtained by observing that if a compact group $G$ has a finitely generated infinite image then either $G/G^0$ of $G^0$ has such image, where $G^0$ denotes the connected component of the identity in $G$.
Now $G/G^0$ is a profinite group, while $G^0$ is a pro-Lie group \cite{Hoffmann}; in the latter case an analogue of Proposition \ref{ss} below for Lie groups suffices to complete the proof.

So let us assume that $G$ is a profinite group with infinite finitely generated image $I$. By considering the closed subgroup of $G$ containing preimages of the generators  of $I$ we reduce to the case when $G$ is topologically finitely generated. Let $I_1$ be the intersection of all subgroups of finite index in $I$. Then $I/I_1$ is a residually finite countable image of $G$ and therefore finite by the argument at the beginning of this section. Thus by replacing $I$ with $I_1$ and $G$ with the preimage of $I_1$ we may further assume that $I$ has no finite images.

In order to apply Corollary \ref{im} we need to be able to understand the countable images of $G/G_0$, which is a semisimple by soluble profinite group.

\begin{proposition}
\label{ss} Let $G=\prod_{i\in X}S_{i}$ be a finitely generated semisimple
profinite group.
Then any infinite image of $G$ is uncountable.
\end{proposition}
Assuming this we can complete the proof of Theorem \ref{images}. Let $N$ be a normal subgroup of $G$ such that $G/N \simeq I$ is a finitely generated countable group without finite images. Then $I=I'$ and hence $G'N=G$. Suppose that $NG_0 <G$, then $G/G_0$ maps onto a nontrivial factors of $I$ which must be infinite.  It follows that the group $\bar G=G/G_0$ has an infinite finitely generated perfect quotient. Let $V$ be the semisimple normal subgroup of $\bar G$ such that $\bar G/V$ is soluble of derived length 3. Clearly $G/V$ cannot have a nontrivial perfect quotient, therefore $V$ must have an infinite finitely generated quotient. This contradicts Proposition \ref{ss} and Theorem \ref{images} follows.

Finally let us indicate a proof of Proposition \ref{ss}.
This proceeds by describing the maximal normal subgroups of a semisimple profinite group $G=\prod_{i\in X}S_{i}$. For this we
need the following result of Liebeck and Shalev:

\begin{theorem}[\cite{ls}]
\label{LS} There is a constant $c$ with the following property: If $S$ is a
nonabelian finite simple group and $C$ is a nontrivial conjugacy class of $S$
then $C^{n}=S$ for some integer $n \leq c \log|S|/\log|C|$.
\end{theorem}

Now let $\mathcal{U}$ be an ultrafilter on the index set $X$. We define a
normal subgroup $N_{\mathcal{U}}$ of $G$ as follows: first define a function
$h:G\rightarrow\lbrack0,1]$ by
\[
h((g_{i}))=\lim_{\mathcal{U}}\frac{\log|g_{i}^{S_{i}}|}{\log|S_{i}|}%
\]
where $\lim_{\mathcal{U}}$ denotes the ultralimit with respect to
$\mathcal{U}$, then set $N_{\mathcal{U}}=h_{\mathcal{U}}^{-1}(0)$. Using
Theorem \ref{LS}, one shows that the subgroups $N_{\mathcal{U}}$ are precisely
the maximal proper normal subgroups of $G$.
Proposition \ref{ss} is then deduced by eliminating the possibility $H\leq
N_{\mathcal{U}}$ when $G/H$ is countably infinite.

\section{An application: locally compact totally disconnected groups}

Profinite groups appear as open subgroup of locally compact totally disconnected (abbreviated lctd) groups. By contrast with profinite groups Corollary \ref{cor1} does not hold even for very 'nice' lctd groups as the following example noticed by D. Segal shows. Let $G$ be the group
\[ G= \left\{ \left( \begin{array}{ccc} 1 & a & c \\ 0 & 1 & b \\ 0 & 0 & 1 \end{array} \right) \ \left| \right. \  a, b \in \mathbb Z, c \in \mathbb Z_p \right\} \] where $\mathbb Z$ has the discrete topology while $\mathbb Z_p$ has the $p$-adic topology. Then $G$ and its compact open subgroup $\mathbb Z_p$ are both finitely generated topologically, however $G'=\mathbb Z < \mathbb Z_p$ is not closed in $G$. 

However the situation may be different for topologically simple tdlc groups (i.e. groups which do not have a nontrivial normal \emph{closed} subgroups). Perhaps the most natural question to ask is whether these groups are abstractly simple.
Note that any non-trivial normal subgroup of a topologically simple group must be dense. In \cite{capmon} Caprace and Monod pose the following
\begin{question} Is there a compactly generated, topologically simple locally compact group which has a proper dense normal subgroup? 
\end{question}
The methods from the previous section provide some information in the case when the topologically simple group is compactly generated and locally finitely generated (i.e its open compact subgroups are finitely generated).

\begin{proposition}\label{tdlc} Let $L$ be a compactly generated, locally finitely generated tdlc group without compact normal subgroups. Let $G$ be a compact open subgroup of $L$ and let $N$ be a dense normal subgroup of $L$. Then $N \geq G'$. In particular $L$ is abstracty simple if and only if $L$ is abstractly perfect.
\end{proposition}
\textbf{Proof:} Suppose that $L$ is generated by a compact set $K$. Without loss of generality we may assume that $K=K^{-1}=GKG$. Let $\mathcal G$ be the Schreier graph of $L$ with respect to $K$ and $U$. The vertices of $\mathcal G$ are $\{Ga\ |\ a \in L\}$  and the edges are the pairs $(Ga,Gka)$ for $k \in K$. Our group $L$ acts transitively on the vertices of $\mathcal G$ with stablizer of the vertex $v_0=G$ equal to $G$. By assumption $L$ has no non-trivial compact normal subgroups hence the action of $L$ on $G$ is faithful. The valency at every vertex of $\mathcal G$ is $|K:U|< \infty$, in particular there are only finitely possibilities for the upper composition factors of $G$ (that is the composition factors of $G/N$ for open normal subgroups $N$ of $G$). We conclude that $G/G_0$ is finite while $G/G'$ has only finitely many nontrivial Sylow subgroups. It follows that there is a finite index subgroup $G>F>G'$ such that for a closed subgroup $M$ of $G$ $G=MF$ implies $G=MG'$.   

Suppose $N$ is a dense normal subgroup of $L$ and put $H=G \cap N$, then $H$ is a dense normal subgroup of $G$. Since $F$ and $G_0$ are open subgroups of $G$ we have that $G=G_0H=FH$. Thus we can find finitely many elements $y_1, \ldots y_r \in H$ such that 
\[ \overline{ \langle y_1, \ldots y_r \rangle}G_0=G=\overline{ \langle y_1, \ldots y_r \rangle}F \] and the last equality implies that
$\overline{ \langle y_1, \ldots y_r \rangle}G'=G$. We can now deduce as in the proof of Corollary \ref{im} that $H>[G_0,G]$. Finally
$G'=[G,HG_0] \leq H[G,G_0] \leq H$ as before. Since $N$ is dense in $L$ we have that $L/N= GN/N=G/H$ is abelian. $\square$ \medskip

In \cite{barnea} the authors study commensurators of profinite groups and one aspect of this is the following question: Which profinite groups can be the compact open subgroups of topologically simple tdlc groups? By the proof of Proposition \ref{tdlc} such a profinite group $G$ has only finitely many different upper composition factors. In particular $G$ cannot be a free profinite group. It is shown in \cite{barnea} that the above question has a positive answer for the Grigorchuk group and negative for the Nottingham group. Another negative answer was proved by G. Willis \cite{willis}: a soluble profinite group cannot be an open subgroup of a compactly generated simple tdlc group, moreover the condition of being compactly generated is necessary. The following concrete question is asked in \cite{barnea}.
\begin{question} Let $G$ be the free pro-$p$ on 2 generators. Can $G$ be an open subgroup in a topologically simple lctd group? 
\end{question}

Next we discuss some of the main ingredients in the proof of Theorem \ref{fin}.

\section{Elements of the proof of Theorem \ref{fin}}
The most natural approach is to use induction on $|\Gamma|$. Suppose for example that
$H$ in Theorem \ref{fin} is an minimal normal subgroup of $\Gamma$, with the property that $H=[H,\Gamma]$. 

Then a necessary
condition on the $y_{i}$ is that at least one of the sets $[H,y_{i}]$ is not
too small, i.e. $y_{i}$ does not centralize a subgroup of $H$ of big Hausdorff
dimension. Now either $H$ is elementary abelian and then each $y_{i}$ acts on
$H$ as a linear transformation, or $H$ is a direct product of isomorphic
simple groups and $y_{i}$ acts by permuting the factors and twisting them with
some automorphisms. This motivates the following definition:

Let the group $\langle Y\rangle$ generated by a set $Y$ act on either \medskip

(1) a set $\Omega$, or

(2) a vector space $V$.$\medskip$

\noindent In situation (1) we say that $Y$ has $(\epsilon,k)$ \emph{fixed
point property} on $\Omega$ if at least $k$ of the elements of $Y$ fix at most
$(1-\epsilon)|\Omega|$ of the points in $\Omega$. Similarly in situation (2)
$Y$ has the $(\epsilon,k)$ \emph{fixed space property} on $V$ if at least $k$
of the elements of $Y$ have centralizer of dimension at most $(1-\epsilon)\dim
V$.

One of the main new ingredients is a result which guarantees that a set
$Y\subset \Gamma$ has the $(\epsilon,k)$-fixed point and fixed space property on
the non-central chief factors of $\Gamma$ provided $Y$ generates $\Gamma$ modulo $\Gamma^{\prime}$ and
modulo $\Gamma_{0}$:

\begin{theorem}
\label{fppThm} Let $\Gamma$ be a finite group with a subset $Y$ such that $\Gamma=\Gamma^{\prime}\left\langle Y\right\rangle =\Gamma_{0}
\left\langle Y\right\rangle $. Then $Y$ has the $\varepsilon/2$-fsp on every
non-central abelian chief factor of $\Gamma$ and the $\varepsilon$-fgp on every
non-abelian chief factor of $\Gamma$ inside $\Gamma_{0}$, where%
\[
\varepsilon=\min\left\{  \frac{1}{1+6\delta},\frac{1}{|Y|}\right\}  .
\]

\end{theorem}

This allows to prove the following.

\begin{theorem}
\label{ThmA}Let $G$ be a group and $K\leq G_{0}$ a normal subgroup of $G$.
Suppose that $G=K\left\langle y_{1},,\ldots,y_{r}\right\rangle =G^{\prime
}\left\langle y_{1},,\ldots,y_{r}\right\rangle $. Then there exist elements
$x_{ij}\in K$ such that%
\[
G=\left\langle y_{i}^{x_{ij}}\mid i=1,\ldots,r,~j=1,\ldots,f_{0}\right\rangle
\]
where $f_{0}=f_{0}(r,\mathrm{d}(G))=O(r\mathrm{d}(G)^{2})$.
\end{theorem}

Let us now consider the situtation where $H$ in Theorem \ref{fin} has a minimal elementary abelian subgroup $M \vartriangleleft \Gamma$ with $[M,\Gamma]=M$. Suppose $h \in H$. 
By induction and considering the smaller group $\Gamma/M$ we may assume that we have found elements $a_{i,j} \in H$ such that for some $m \in M$ we have
\[ hm= \prod_{i=1}^f \prod_{j=1}^r [a_{i,j},y_j] \]
Let us replace $a_{i,j}$ by $a_{i,j}x_{i,j}$ with some $x_{i,j} \in M$. By collecting the terms involving $a_{i,j}$ to the left we reach the equation
\begin{equation} \label{ee} m= \prod_{i=1}^f \prod_{j=1}^r [x'_{i,j},z_{i,j}] \end{equation} where each $z_{i,j}$ is specific conjugates of $y_j$ depending on $x_{i,j}$ and $a_{i,j}$ and similarly for $x'_{i,j}$.
In order to be able to solve this equation in terms of $x'_{i,j}$ we assume that the elements $y_{i,j}$ generate $G/M$ and put this condition in the induction hypothesis for $\Gamma$, see Key Theorem in \cite{NS} Theorem 3.10. But then when solving (\ref{ee}) for $M$ we need to ensure that the extra condition that $z_{i,j}$ generate $\Gamma$ is satisfiied. We do this by counting the solutions of (\ref{ee}) and comparing this with the number of conjuagtes of $y_{i}$ whcih fail to generated $\Gamma$.
The counting is relatively straightforward when $M$ is a minimal normal subgroup as above (i.e. $M$ is abelian and $M=[M, \Gamma]$) but if $[M,\Gamma]=1$ or when $M$ is semisimple it gets much more complicated. To deal with the situation when $M \leq Z(\Gamma)$ the methods developed earlier in \cite{NS3} suffice. When $M$ is semisimple we need better lower estimates for the number of solutions to the equation (\ref{ee}) than the ones in \cite{NS3}.
\medskip

This leads us to the subject of growth in finite simple groups and the following generalization of the `Gowers
trick' by Babai, Nikolov and Pyber in \cite{BNP}:

\begin{theorem}
For a finite group $G$ let $l=l(G)$ denote the minimal degree of a nontrivial
real representation of $G$. Suppose now that $k \geq3$ and $A_{1}, \ldots,
A_{k}$ are subsets of $G$ with the property that $\prod_{i=1}^{k} |A_{i}|
\geq|G|^{k}/l^{k-2}$. Then $\prod_{i=1}^{k} A_{i}=G$.
\end{theorem}

This is a key ingredient in the proof of Theorem \ref{comm20} below. To make
sure that the conditions of the Gowers trick are met we prove the following
result concerning twisted commutators of finite simple groups:

For automorphisms $\alpha,\beta$ of a group $G$ and $x,y\in S$ we write%
\[
T_{\alpha,\beta}(x,y)=x^{-1}y^{-1}x^{\alpha}y^{\beta}.
\]

\begin{theorem}
\label{T1}There exist $\varepsilon>0$ and $D\in\mathbb{N}$ such that if $S$ is
a finite quasisimple group with $l=l(S)>2$, $\alpha,\beta\in\mathrm{Aut}%
(S)^{(D)}$, and $X\subseteq S^{(2D)}$ has size at least $(1-\varepsilon
)\left\vert S^{(2D)}\right\vert $, then
\[
\left\vert \mathbf{T}_{\alpha,\beta}(X)\right\vert \geq l^{-3/5} |S|
\]

\end{theorem}

Note that there are only finitely many quasisimple groups with $l(S)\leq2$.

The last two theorems combined with some further work give

\begin{theorem}
\label{comm20}Let $D$ and $\epsilon$ be the constants introduced in Theorem
\ref{T1}. Let $N$ be a finite quasisemisimple group with at least $3$
non-abelian composition factors. Let $\mathbf{y}_{1},\ldots,\mathbf{y}_{10}$
be $m$-tuples of automorphisms of $N$. Assume that for each $i$, the group
$\left\langle \mathbf{y}_{i}\right\rangle $ permutes the set $\Omega$ of
quasisimple factors of $N$ transitively and that $\mathbf{y}_{i}$ has the
$(k,\eta)$-fpp on $\Omega$, where $k\eta\geq4+2D$. For each \thinspace$i$ let
$W(i)\subseteq N^{(m)}$ be a subset with $\left\vert W(i)\right\vert
\geq(1-\varepsilon/6)\left\vert N\right\vert ^{m}$. Then%
\[%
{\displaystyle\prod\limits_{i=1}^{10}}
W(i)\phi(i)=N
\]
where $\phi(i):N^{(m)}\rightarrow N$ is given by%
\[
(x_{1},\ldots,x_{m})\phi(i)=%
{\displaystyle\prod\limits_{i=1}^{m}}
[x_{i},y_{ij}].
\]

\end{theorem}

This last theorem provides the induction step in the
proof of Theorem \ref{fin} in the case when $M\leq H$ is a nonabelian minimal normal
subgroup of $G$.
\section{Growth in finite simple groups}
A common theme in this survey has been product decompositions of finite and profinite groups, and how these relate to algebraic properties of infinite groups.  
Most useful are results of the following kind: 

Suppose that $\Gamma$ is a finite group and $X_1, \ldots X_n$ are subsets of $\Gamma$. If 
\[ \sum_{i=1}^n \log |X_i|>C \log|\Gamma|\] for some constant $C$ and $X_i$ generate $G$ we want to deduce that $X_1 \cdots X_n=\Gamma$.  
Of course this result in not true in all finite groups, the cyclic abelian groups being a counterexample.

However when we turn to the other extreme, finite simple groups then we expect positive answers.
For example if $X_i$ are all equal to a conjugacy class in $\Gamma$ and $\Gamma$ is a finite simple group this is the content of Theorem \ref{LS}.

A very influential problem in the area is the Babai Conjecture \cite{BS}:
\begin{conjecture} There is a constant $C>0$ such that if $S$ is a finite simple group and $X=X^{-1}$ is a generating set for $S$ then the diameter of the Cayley graph $Cay(S,X)$ is at most $(\log |S|)^C$.
\end{conjecture}
Note that when $1 \in X$ the Babai conjecture is equivalent to $S=X^{*n}$ for some integer $n < (\log |S|)^C$.

The conjecture was proved in the positive by Helfgott \cite{Helfgott} for $PSL_2(\mathbb F_p)$ and subsequently extented in \cite{Helfgott2} to $PSL_3$ and other non-prime fields by Dinai \cite{dinai}. Recently Breuillard, Green and Tao \cite{BGT} and at the same time Pyber and Szabo \cite{PS} have proved the Babai Conjecture for all finite simple group of bounded Lie rank. Both papers prove a \emph{Helfgott type} estimate on the growth of subsets in simple groups. Rather than define what we mean by this we give a statement of one of the main theorem in \cite{PS}, or equivalently Corollary 2.4 in \cite{BGT}.
\begin{theorem} Let $S$ be a finite simple group of Lie type of rank $r$ and
$A$ a generating set of $L$. Then either $AAA = L$ or
$|AAA| > \delta|A|^{1+\epsilon}$
where $\epsilon, \delta$ depend only on $r$.
\end{theorem}

Results of this type together with sieve methods developed by Bourgain, Gamburd and Sarnak \cite{BGS} play an important role in the construction of expanders as Cayley graphs in finite simple groups of bounded Lie rank. For example this was how the Suzuki finite simple groups were shown to be a family of expanders in \cite{BGT1}. For details about this fascinating and rapidly evolving subject we point the reader to the survey by B. Green \cite{green}. 

Expressing a finite simple group as a product of few of its subgroups has been used to construct expanders in large rank as well, see \cite{KLN} and the references therein. Let us mention a very general conjecture made in \cite{LNS}.
\begin{conjecture} There exists an absolute constant $c$ such that if $S$ is a finite
simple group and $A$ is any subset of $S$ of size at least two, then $S$ is a
product of $N$ conjugates of $A$ for some $N < c \log |S|/ \log |A|$.
\end{conjecture}

This conjecture can be viewed as a direct generalization of Theorem \ref{LS} from conjugacy classes to subsets. Some special cases are known, for example it holds for specific subgroups $A$ of $G$, or when $|A|$ is bounded. Most notably its validity has recently been proved for simple groups of bounded rank in \cite{Pybermany}.  

\section{Rank gradient}

We saw that some properties of $\Gamma$ cannot be deduced from its profinite completion $\hat \Gamma$, for example $\hat \Gamma$ does not determine $\Gamma$ up to an isomorphism. Another such example is $d(\Gamma)$, the minimal number of generators of $\Gamma$, which cannot be found from knowledge of $d (\hat \Gamma)$ alone: For any integer $n$ there exists a residually finite group $\Gamma$ such that $d(\Gamma)=n$ but $\hat \Gamma$ (equivalently every finite image of $\Gamma$) is $3$-generated, see \cite{Wise}. It is thus surprising that the growth of $d(H)$ for subgroups of finite index in $\Gamma$ can be recovered from the knowledge of $\hat \Gamma$ plus some extra information: the action of $\Gamma$ on its profinite completion. Let us make this precise. \medskip

Let $\Gamma$ be a finitely generated group and $(\Gamma_i)$ a chain of subgroups in $\Gamma$.
The \emph{rank gradient} of $\Gamma$ with respect to $(\Gamma_i)$ is defined as
 \[ \mathrm{RG}(\Gamma,(\Gamma_i))= \lim_{i \rightarrow \infty} \frac{d(\Gamma_i)-1}{[\Gamma:\Gamma_i]}\]
 where $d(\Gamma)$ denotes the minimal number of generators of $\Gamma$.
This notion has been introduced by Marc Lackenby \cite{lack}
in the study of hyperbolic 3-manifolds and the virtually Haken conjecture.
A natural question is whether the rank gradient depends on the choice of
the normal chain in $\Gamma$.

\begin{conjecture}\label{rg} If $(\Gamma_i)$ and $(\Delta_i)$ are two normal chains in $\Gamma$ with
trivial intersection then $\mathrm{RG}(\Gamma, (\Gamma_i))= \mathrm{RG}(\Gamma, (\Delta_i))$.
\end{conjecture} 

Conjecture \ref{rg} is relevant the following well-known problem in
3-dimensional topology which dates back to Waldhausen: \emph{Is there a closed hyperbolic $3$-manifold $M$ such that $d(\pi_1(M)) \not = g(M)$? (where $\pi_1(M)$ is the fundamental group of $M$ and $g(M)$ is the Heegaard genus of $M$)} 
This has recently been solved by T. Li:

\begin{theorem} [\cite{Li}] For any fiven $n \in \mathbb N$ there exists a closed orientable hyperbolic 3-manifold $M$ such that $g(M)- d(\pi_1(M)) >n.$
\end{theorem}

As explained in \cite{AN} the truth of Conjecture \ref{rg} together
with results of M. Lackenby \cite{lack1} and A. Reid \cite{AR} implies an even stronger result: 
the ratio $d(\pi_1(M))/g(M)$ can be arbitrarily small for closed hyperbolic
3-manifolds $M$. For this application it will be enough to prove Conjecture 1 for free-by-cyclic groups.

Conjecture \ref{rg} is closely related to the subject of \textrm{measurable group actions} and in particular the notion of \emph{cost} as introduced by Levitt and developed by Gaboriau in \cite{gabor}.
Let $X$ be a probability measure space. Let $\Gamma$ be a group acting on
$X$ by measure preserving transformations. Assume that the action of $\Gamma$
on $X$ is ergodic and essentially free. As explained in \cite{gabor} and also in \cite{KMiller} one can define an invariant $\mathrm{cost}(\Gamma,X)$, the cost of the action of $\Gamma$ on $X$. This invariant has been used by Gaboriau to prove that
two finitely generated free groups have measure equivalent actions if and only if they are isomorphic. It turns out that rank gradient is connected to cost via the following result from \cite{AN}:

\begin{theorem}\label{cost} Suppose $(\Gamma_i)$ is a normal chain with trivial
intersection in a finitely generated group $\Gamma$. Let 
$\widehat \Gamma = \underleftarrow{\lim}_i \Gamma/\Gamma_i$ be the completion of $\Gamma$ with respect
to $(\Gamma_i)$ and let $\Gamma$ act on $\widehat \Gamma$ by left translations. Then
\[ \mathrm{RG}(\Gamma,(\Gamma_i)) = \mathrm{cost}(\Gamma, \widehat \Gamma) -1. \]  
\end{theorem} 

There is not a single group known which has two (essentially free) measurable actions
with different cost. This is the content of the
\medskip

\noindent \textbf{Fixed Price Conjecture} \emph{Every countable group has the same
cost in each of its measurable essentially free actions on probability spaces.}

It is known \cite{gabor} (see also \cite{KMiller}) that the collection of groups with fixed price,
includes free groups, amenable groups, groups with infinite centre and is closed under free products.
In view of Theorem \ref{cost} the  validity of the Fixed Price Conjecture will imply Conjecture \ref{rg}.  \medskip

\emph{Acknowledgements.} Thanks to D. Segal from whom the author has learned much about profinite groups and who made many helpful comments on a preliminary version of this article. Y. Barnea and A. Shalev provided up to date information for Sections 5 and 7. \bigskip 

\bigskip

Department of Mathematics,

Imperial College London, 

SW7 2AZ, UK. 

n.nikolov@imperial.ac.uk

\end{document}